\documentclass[reqno,12pt]{amsart}
\usepackage[indentafter]{titlesec}
\titleformat{name=\section}{}{\thetitle.}{0.8em}{\bfseries}
\titleformat{name=\subsection}[runin]{}{\thetitle.}{0.5em}{\bfseries}[.]
\titleformat{name=\subsubsection}[runin]{}{\thetitle.}{0.5em}{\bfseries}[.]
\usepackage{hyperref}
\usepackage{cleveref}
\usepackage[english]{babel}
\usepackage[all,cmtip]{xy}
\usepackage{lmodern}
\newtheorem{thm}{Theorem}[section]
 \numberwithin{dummy}{section}

\theoremstyle{definition}
\usepackage{graphicx,color}
\usepackage{graphics}
\usepackage{amssymb}
\usepackage{subfloat}
\usepackage{float}
\usepackage{array,multirow,makecell}
\usepackage{tabularx}
\usepackage[final]{pdfpages}
\usepackage{amsmath,amsfonts,amssymb}
\usepackage{systeme}
\usepackage{cases}
\usepackage{algorithmic}
\usepackage[]{algorithm2e}
 \usepackage{easybmat} 
\usepackage{fancyvrb}
\usepackage{verbatim}
\usepackage{xcolor}
\usepackage{caption}
\usepackage{mathpazo} 
 
 \theoremstyle{remark}
 \newtheorem{rem}[thm]{\textbf{Remark}}
\theoremstyle{plain}

\newtheorem{theorem}{Theorem}[section]

\newtheorem*{pof*}{Proof of Lemma}

\newtheorem*{pofd*}{Proof}
\newtheorem*{pofc*}{Proof of corrolary}
\newtheorem*{proof of theorem}{Proof of Theorem}
\newtheorem{e-proposition}[theorem]{Proposition}
\newtheorem{corollary}[theorem]{Corollary}
 
\usepackage[left=3cm, right=3cm, bottom=3cm]{geometry}
\usepackage{fancyhdr}
\newcommand\raisepunct[1]{\,\mathpunct{\raisebox{0.5ex}{#1}}}

\author{R.Malek$^{(*)}$ \& C.Ziti $^{(\P)}$}
\address{$^{(*)}$ Moulay Ismail University of Meknes, Faculty of Sciences. \vskip 3pt
	Team work: Equations aux Dérivées Partielles et Calcul Scientifique, (EDPCS).}
\email{\textcolor[rgb]{0.00,0.00,1.00}{r.malek@edu.umi.ac.ma}}

\address{$^{(\P)}$ Moulay Ismail University of Meknes, Faculty of Sciences. \vskip 3pt
Team work: Equations aux Dérivées Partielles et Calcul Scientifique, (EDPCS).}
\email{\textcolor[rgb]{0.00,0.00,1.00}{chziti@gmail.com}}

\title{ Extension of $\delta_{-}ziti$ method in the unit ball:\\Numerical integration, resolution of Poisson's problem and  Heat transfer.}

\begin{document}

\medskip
\medskip
\medskip
\maketitle

\vspace*{-5mm}

%
%
\bigskip
\begin{abstract} 
Inspired by the Galerkin and  particular method, a new approximation approach is recalled in the Cartesian case. In this paper, we are interested specially by constructing this method, when the domain of consideration is a two dimensional ball, to extend this work to the several dimension.  We reduce the number of iterations to calculate integrals and numerical solution of Poisson and the Heat problems (elliptic nd parabolic PDEs), in a very fast way.\\
\end{abstract}

\section{\textbf{ Introduction} } \par 

The chemotactic dynamics of a population requires several steps, particularly,  aggregation and blow-up. The Keller-Segel model describes this phenomena. It was suggested  by Patlak in 1953 \cite{B}, Keller-Segel in 1970 \cite{C}, which allows for both diffusion and aggregation: depending on the initial data, the solution might exist globally in time or blow up in finite time, depending on the balance of forces between the different parameters involved in the system, the blow-up phenomenon may or may not occur. \vskip 4pt 
In fact, the blow-up is a singular behaviour of a Dirac solution. Most of numerical methods ( e.g.\, Galerkin, Particular method, spectral method..) does not ensure the transition from a regular behaviour to another singular one (\textit{i.e.}\, detection of blow-up ). Under certain formulations of the Keller-Segel  model, the phenomenon of aggregation has been shown to lead to finite-time blow-up. A large body of works has been devoted to determine when blow-up occurs or whether globally solutions exists:  Authors of \cite{D} developed a family of new interior penalty discontinuous Galerkin methods for solving the Keller–Segel chemotaxis model. In  \cite{E} they  investigated non-negativity of exact and numerical solutions to a generalized Keller–Segel model where this model includes the so-called minimal Keller–Segel model.  The main aim of \cite{F} is to study the Keller–Segel model of chemotaxis and to develop a composite particle-grid numerical method with adaptive time stepping which allows to resolve and propagate singular solutions. 
The purpose of \cite{G} is to formulate a phenomenological model from which the existence and properties of migrating bands can be deduced. Authors of \cite{H} detected the blow-up as $\delta$-function (amoebae aggregation) at the proximity of the origin in dimension one and  on a ball in a multidimensional space. Therefore, it was necessary to find a new numerical scheme, which detect this type of singularities easily, without loosing the advantages of classical methods. The $\delta_{-} ziti$ method is on the challenge. It was tested on several type of problems (see \cite{Aa} and \cite{A}), including the Keller-Segel model, but only in the Cartesian case ( segments, rectangle, cube $\cdots$ ).\vskip 6pt  
 The main goal of  $\delta_{-} ziti$ method, is to approach a function with several variables,  to integrate it in a given domain, and  to resolve numerically Partial and Ordinate Differential Equations (PDEs and ODEs). This method is based on the classical variation formulation of Galerkin and  the most important step, is the construction of our orthonormal family, from the famous function  $\Phi$  $\in$ $D(\Omega)$ with compact support, defined by  
\begin{equation}
\Phi(x)=
\begin{cases}
\exp(\frac{1}{\mid x\mid^2 -R^2}) & \text{si $\mid x\mid < R$} \\
0  & \text{otherwise,}
\end{cases} 
\label{Phi}
\end{equation}
where R $ > $ 0, $\Omega \subset \mathbb{R}^n$ and $x \in \mathbb{R}^n$.\\
This function is used especially in numerical analysis, distributions and  functional analysis. It is characterized by giving the best approximation of the Dirac measure.
In \cite{Aa} and \cite{A}, the multi-dimensional Cartesian case was detailed. \vskip 6pt
The main aim of this paper is the construction of $\delta_{-} ziti$ method when the domain is a disk in the two-dimensional case, (in general, a multi-dimensional ball). To generalize this method, we opt for two strategies: the first one consists in sweeping all the disk with segments, in the two directions, as shown in figure \ref{fig1} and  to reconstruct our basis functions in every segment, which means that we inject all the work already done in the mono-dimensional case. To test this strategy, we apply the resulting tools to calculate numerically integrals and to solve partial differential equations (two tests will be detailed; an elliptic equation "The Poisson problem" and a parabolic one "The Heat equation"). The second strategy is a direct use of the polar parametrisation of a disk, we will show that this strategy is also efficient and gives us a good approximation ( Integrals ans two tests of resolving PDEs). \vskip 6pt

The outline of this paper is as follows. In section \ref{cartt}, we will present the mathematical tools of construction, which permits to apply this method in the Cartesian case, as shown in figure  \ref{fig1},  to calculate, numerically, some integrals defined in a disk domain.\vskip 6pt
The section \ref{poll} is devoted to the construction of the method's fundamental elements, using polar coordinates. Like the previous section, one of the most important parts is the numerical integration using our method and  in the two cases we will compare the exact value of an integral, by the numerical one, obtained by $\delta_{-} ziti$ method.\vskip 6pt
In the last section \ref{EDP}, we apply our approach to find the numerical solution of Poisson problem and  the Heat equation. Our goal is to compare the solution obtained by $\delta_{-} ziti$ method, with a given analytical one defined in a disk domain and  to calculate the error in $L^\infty(\Omega)$. By the next, we present an approximated solution using the finite element method and  we compare it with our one.
\section{\textbf{Overview of the mono-dimensional construction.  } }  \vskip 6pt
All the results presents in this section, are proved in \cite{Aa} and \cite{A}. The  fundamental results of construction are given as follow:\vskip 6pt
First, we take a uniform mesh of [a,b]  with the step $h=\frac{b-a}{N}$, where N is
an integer such that
$x_i=a+(i-1)h,\ \forall i \in \left[1,N+1\right].$\\
From the function $\Phi$ defined in (\ref{Phi}), we define  $\varphi_\epsilon$ by:
$$ \varphi_\epsilon(x)=\displaystyle \frac{c}{\epsilon}\Phi\left(\frac{x}{\epsilon}\right), \ \ \ \ \text{for all}~~ \epsilon > 0,$$
where c:=$\frac{1}{\int_{\mathbb{R}} \Phi(x)dx}$ is the constant of normalization. \par 
This sequence  $\varphi_\epsilon$ converges to Dirac in the sense of distributions, which is often used to detect singularities.\par 
We construct the family  $(\varphi_i)_{i=1 \cdots N+1}$ as follows:
\begin{equation*}
\left\lbrace
\begin{aligned} \label{phi}
&\varphi_i(x)=\varphi_h(x-x_i)=\frac{C}{h}\Phi(\frac{x-x_i}{h}),~~~~ \text{for all} \ \  x\in[x_{i-1},x_{i+1}], \ \  \  i \in[2,N], \\
&\varphi_1(x)=\varphi_h(x-x_1)=\frac{C}{h}\Phi(\frac{x-x_1}{h}),~~~~ \text{for all} \ \  x\in[x_{1},x_{2}],\\
&\varphi_{N+1}(x)=\varphi_h(x-x_{N+1})=\frac{C}{h}\Phi(\frac{x-x_{N+1}}{h}), ~~~~\text{for all} \ \  x\in[x_{N},x_{N+1}].
\end{aligned}
\right. 
\end{equation*}
Let consider the Hilbert space $L^2(\mathbb{R})$, with the usual scalar product 
$(~ ~, ~ ~ )$. Observe that the family $\left(\varphi_i \right)_{1\leq i \leq N+1}$ is linearly independent, then using the Gram-Schmidt process, we construct a unique orthogonal family, noted $ \left( \tilde{\Psi }_i \right)$ satisfying the following relation

\begin{equation}
\left\lbrace
\begin{aligned} \label{psi}
& \tilde{\Psi}_{i}(x) = \varphi_{i}(x)+\lambda_{i-1} \tilde{\Psi}_{i-1}(x),\\
& \lambda_{1}=-\frac{\alpha}{\beta},\\
& \lambda_{i+1}=g(\lambda_{i}),\\
\end{aligned}
\right. 
\end{equation}
with $$g(X)=\frac{\lambda_{1}}{2-\lambda_{1} X}, \ \ \alpha=(\varphi_{1},\varphi_{2}), \ \ \beta=(\varphi_{1},\varphi_{1}).$$
The spectral method applied to find the direct formula of the basis functions gives, 
\begin{equation*}
\lambda_{i}=-\frac{(\varphi_{i},\psi_{i-1})}{(\psi_{i},\psi_{i-1})}\raisepunct{,}
\end{equation*}
and the recurrence application of the definition given in (\ref{psi}) gives the following formula:
\begin{equation*}
\tilde{\Psi}_{i}(x) =\varphi_{i}(x)+\lambda_{i-1}\varphi_{i-1}(x)+\lambda_{i-1}\lambda_{i-2}\varphi_{i-2}(x)+\cdots+\prod_{k=i-1}^{1} \lambda_{k} \varphi_{1}.
\end{equation*}
Let  $\Psi_i= \frac{\tilde{\Psi}_{i}}{|| \tilde{\Psi_i}||}$ the normalization of  $\tilde{\Psi}_{i}$ ( for more details see \cite{Aa} and \cite{A}).\\
The method permits to approach a given function $f$ and  an integral, using the following relations:
\begin{align}
\begin{split}
& f(x)  \simeq \sum_{i=1}^{N} c_i \Psi_{i}(x),\\
& c_i \simeq  \int_{a}^{b} f(x)  \Psi_{i}(x)dx,\\
& \int_{a}^{b} f(x) dx \simeq \sum_{i=1}^{N} c_i I_i,
\label{12}
\end{split}
\end{align}
where $ I_i:=\int_{a}^{b}\Psi_{i}(x)dx$. If we take $x=r_k$ in (\ref{12}), we obtain:
\begin{eqnarray}
\begin{split}
& c_i \simeq  \displaystyle \frac{f(r_i)}{\Psi_i(r_i)} \raisepunct{,} i=1\cdots N-1,\\
& c_{N} \simeq \displaystyle \frac{f(b)}{\Psi_N(b)}\raisepunct{,}\\
&  \int_{a}^{b} f(x)  \Psi_{i}(x)dx\simeq \displaystyle \frac{f(r_i)}{\Psi_i(r_i)}\raisepunct{,}\\
& \int_{a}^{b} f(x) dx \simeq \sum_{i=1}^{N-1} \displaystyle \frac{f(r_i)}{\Psi_i^2(r_i)}+ \frac{f(b)}{\Psi_{N}(b)^2}\raisepunct{.}
\label{34}
\end{split}
\end{eqnarray}

\parindent=0em 
To reduce the iterations number, in \cite{these}  they proved that,
\begin{equation*}
|\lambda_{i+1}-\lambda_{i}|<\epsilon \ \ \text{as soon as }\ \  i\geq N_0=\displaystyle \left[\frac{\ln(\displaystyle \frac{ \epsilon(2-\lambda_{1}^2)}{\lambda_{1}^3-\lambda_{1}})}{\ln(\displaystyle\frac{\lambda_{1}}{2+\lambda_{1}})^2}\right]+1,  
\end{equation*}
where $[.]$ denotes the floor function. In particular for  \ \ $\epsilon=10^{-M}$, we concluded that the parameter $\lambda_i$ is nearly stationary from a certain rank, which reduces considerably the number of iterations.
Using $r_i$ as a root of $\Psi_{i+1} $, we can define the parameter $\lambda_{i}$ by, $\lambda_{i}=\displaystyle -\frac{\varphi_{i+1}(r_i)}{\varphi_{i}(r_i)}.$

\section{ \textbf{The first strategy:  Cartesian coordinates.}}\label{cartt}
\subsection{Construction of intern nodes} \leavevmode \par
\begin{figure}[htbp!] 
	\includegraphics[width=13cm]{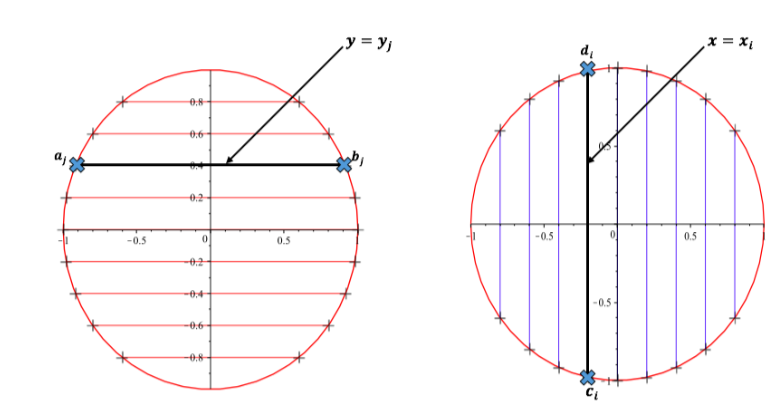}
	\caption*{$a_j=-\sqrt{1-y_j^2},  \ \ b_j=-a_j$ and  $c_i=-\sqrt{1-x_i^2}, \ \ d_i=-c_i$}
	\caption{Horizontal and vertical segments }
	\label{fig1}
	\end{figure}

 
In this section, we are interested in the extension of $\delta_{-} ziti$ method, when $\Omega$ is a disk centred in the point $O=(0,0)$ (or the ball in the multi-dimensional  case). As a first step, we start by a general presentation of this new strategy.\vskip 6pt
We present the important steps of the construction, inspired by the mono-dimensional case. For this, suppose that we can sweep the inside of the domain by a set of intervals, horizontally and vertically, therefore, all the work resides in the construction of the nodes in every interval. (see Figure  \ref{fig1})



An  algorithm which explains the steps of the construction will be presented in the next part . The main idea is to fix the number of nodes in every interval and  to vary the step of subdivision associated to every interval (horizontally and vertically). Every node is noted by $(x_i^j,y_j^i)$, (see  Figure \ref{fig1}). \par
\begin{rem} \leavevmode \par 
	Note that, for every fixed vertical level $j$ (respectively the horizontal level  $i$), every internal segment is limited by $a_j=-\sqrt{1-y_j^2}$ and $b_j=\sqrt{1-y_j^2}$ , (respectively, $c_i=-\sqrt{1-x_i^2}$  and $d_i=\sqrt{1-x_i^2}$ ).\\ 
	For simplicity, the step of horizontal subdivision  will be noted $h_j$  (respectively, the vertical step will be noted $h_i$).
\end{rem} 
Here we present an algorithm to calculate the  internal nodes.\par 
\newpage
\begin{algorithm}
	\KwData \par  {~~~~~~~~~~~~\ \ \ $N=$ The number of nodes, in every interval \\~~~~~~~~~~~~\ \ \  Fix an interval [a,b]  \\ ~~~~~~~~~~~~\ \ \ $h^1=\frac{b-a}{N}$ } \;
	\For {$i=1\cdots N$}{
		$x_{i}^1=a+(i-1)h^1$ horizontal nodes for the first intervall. \\
		$y_{i}^1=a+(i-1)h^1$ vertical nodes for the first intervall
	} 
	
	\For {$j=2\cdots N$}{
		\For {$i=1\cdots N$}{
			a=-$\sqrt{1-(-1+j.h)^2}$\\
			b=-a\\
			$h^j=\frac{b-a}{N}$ \\
			$x_{i}^j=a+(i-1)h^j$ \\
			
		}
	}
 \label{alg1}
 \caption{ Construction of the nodes in the Cartesian case. }
\end{algorithm}
\subsection{Construction of the orthonormal set} \leavevmode \par 
For every node $x^j_i$ (respectively,  $y^i_j$) we associate the function $\varphi_i^j$, ( noted  $\varphi_i$ if there is no ambiguity) (respectively ,the family  $\varphi_j^i$ will be noted  $\varphi_j$) defined by:
\begin{align}
\begin{split}
\varphi_i(x):=\frac{c}{h_j} \Phi(\frac{x-x_i^j}{h^j}), ~~ \forall i,j=1\cdots N, \\
\varphi_j(y):=\frac{c}{h_i} \Phi(\frac{y-y_j^i}{h^i}), ~~ \forall j,i=1\cdots N,
\end{split}
\end{align}
where,\\
$\bullet$ $h^j$ is the step of construction in the horizontal interval of indication $j$, which describe the distance between the node $x_i^j$ and  $x_{i+1}^j$.\\
$\bullet$ $h_i$ is the step of subdivision in the vertical  interval of indication  $i$, which describe the  distance  between  $y_j^i$ and  $y_{j+1}^i$.\\ \\ 
It is simple to see that the family $(\varphi_{i})$ is linearly independent, so we can  construct an orthogonal family $(\tilde{\Psi_i})_{i=1 \cdots N}$ by using the Gram-Schmidt process, in the space $\mathbb{L}^2([a,b])$, (construction  in every internal interval of the domain $\Omega=B(0,1)$, horizontally and vertically), verifying the following relation, 
\begin{equation*} {Horizontally:}
\left\lbrace
\begin{aligned}
\tilde{\Psi}_{1}(x)&=\varphi_1(x) \\         
\tilde{\Psi}_{i}(x)&=\varphi_{i}(x)+ \sum_{k=1}^{i-1}\lambda^{(i)}_k\tilde{\Psi}_k(x) , \ \ \ \  \text{for all} \ \  i=2,\dots,N,\\
\end{aligned}
\right.
\label{defh}
\end{equation*}

\begin{equation*} {Vertically:}
\left\lbrace
\begin{aligned}
\tilde{\Psi}_{1}(y)&=\varphi_1(y) \\         
\tilde{\Psi}_{i}(y)&=\varphi_{i}(y)+ \sum_{k=1}^{i-1}\lambda^{(i)}_k\tilde{\Psi}_k(y) , \ \ \ \  \text{for all} \ \  i=2,\dots,N,\\
\end{aligned}
\right.
\label{defv}
\end{equation*} 

which will be reduced in the following theorem, already proved in the mono-dimensional case, (see \cite{Aa} and \cite{A}).
\begin{theorem}\leavevmode \par 
The orthogonal family $(\tilde{\Psi_i})_{i=1 \cdots N}$ (vertically and horizontally), verify the following recurrence relation:
	\begin{equation}
	\left\lbrace
	\begin{aligned}
	&\tilde{\Psi}_{1}=\varphi_1 \\         
	& \tilde{\Psi}_{i+1}=\varphi_{i+1}+\lambda_i\tilde{\Psi}_i , \ \ \ \  \text{for all} \ \  i=1 \cdots N-1 ,
	\\  
	&\lambda_{i-1}=- \frac{(\varphi_i,\varphi_{i-1})}{(\tilde{\Psi}_{i-1},\tilde{\Psi}_{i-1})}\raisepunct{,}\\
	\end{aligned}
	\right.
	\label{thm2.2}
	\end{equation}
	\label{thm22}
	where $( , )$ is the usual scalar product in the Hilbert space $L^2([a,b])$
\end{theorem}

\begin{corollary} \leavevmode \par 
	The family $(\tilde{\Psi}_{i})$ and the set $(\lambda_{i})$ defined in the theorem ~\eqref{thm22}, verify the following relations:
	\begin{align}
	\begin{split}
& 	1)  \ \ \ \tilde{\Psi}_{i}=\varphi_i+\lambda_{i-1} \varphi_{i-1} +\lambda_{i-1} \lambda_{i-2} \varphi_{i-2}+\cdots+\lambda_{i-1} \cdots \lambda_{1} \varphi_{1}.\\ 
&2) \ \ \  \tilde{\Psi}_{i}(x_i^j)=\varphi_{i}(x_i^j)=\frac{c}{h_j^2e}.\\ 
&3) \ \ \  \tilde{\Psi}_{i}(y_i^j)=\varphi_{i}(y_i^j)=\frac{c}{h_i^2e}.\\ 
&4) \ \ \ \text{In every fixed level}, \  (\varphi_{i},\tilde{\Psi}_{i-1}) = (\varphi_{i},\varphi_{i-1}) .\\
&5)\ \ \   -1<\lambda_{i}=-\frac{(\varphi_{i},\varphi_{i+1})}{(\tilde{\Psi_{i}},\tilde{\Psi_{i}})}<0.
	\end{split}
	\end{align}
	\label{corlm}
\end{corollary}

\subsection{Fundamental results: Numerical integrations}\leavevmode \par 
In this paragraph, we are interested by the approximation of integrals, where the domain  is the unit disk $\Omega=B(0,1)$, using the horizontal and vertical test functions, as well as the roots, verifying the following relations:\par  
\begin{equation}
\begin{split}
& \Psi_{ij}(x,y)=\Psi_i^j(x).\Psi_j^i(y), \ \ \ \forall i,j=1\cdots N,\\
& r_{ij}=(r_i^j,s_j^i),
\end{split}
\end{equation}
where, \\
$\bullet$ $\Psi_{i}^j(x)$ are the basis functions in the horizontal dimension, (respectively, $\Psi_{j}^i(y)$ are the basis functions in the vertical dimension).\\
$\bullet$ $r_i^j$  are the roots of $\Psi_{i}^j(x)$ (respectively, $s_j^i$ are the roots of $\Psi_{j}^i(y)$).\\
In this section, we are interested by the approximation of a double integral, defined in a disk domain, using (\ref{12}), we obtain the following results:
\begin{theorem} \leavevmode \par 
	Let $\Omega=B(0,1)$,  $g$ a given function in $L^2(\Omega)$, and N denotes the roots number, 
	therefore, we have the following approximations:

	\begin{equation}
	\begin{split}
	& \int_{\Omega} g(x,y) \ dx dy \simeq \sum_{i,j=1}^{N} \frac{g(r_i^j,s_j^i)}{\Psi_{i}^j(r_i^j).\Psi_{i}^i(s_j^i)} \int_{a_j}^{b_j} \Psi_{i}^j (x) dx \int_{c_i}^{d_i} \Psi_{j}^i (y) dy, \\
	& \int_{\Omega} g(x,y) \Psi_{ij}(x,y)\  dx dy \simeq  \frac{g(r_i^j,s_j^i)}{(\Psi_{i}^j(r_i^j).\Psi_{i}^i(s_j^i))^2}\raisepunct{.}
	\end{split}
	\end{equation}
	\label{integ}
	here we take, $r_N^j=b_j$ and $s_i^N=b_i$.
\end{theorem} 
In the table (\ref{tab1}), we present some numerical tests of integration.  We compare the exact value, with the numerical approximation, obtained by $\delta_{-} ziti$ method, in the Cartesian case.\\ \\ 
\renewcommand{\arraystretch}{2.5}
\begin{table}[ht]
	\centering
\begin{tabular}{ |c |  c|  c  |   c ||   } 
	\hline \hline 
	$\int_{\Omega} f(x,y) dx dy$ & $2\pi.$Ex& $2\pi.$Ap& $2\pi.$Error \\
	\hline
	$\int_{\Omega} (\frac{1}{x^2+y^2+1})^{\frac{1}{4}} $ &$0.4545285537$ &$0.454496459918650$ & $0.0000320937813496069$  \\ 
	\hline 
	$\int_{\Omega} \sqrt{(\frac{1}{x^2+y^2+1})} $ &$0.4142135624$ &$0.414440692467526$ & $0.0002271300675258380$  \\ 
	\hline 
	$\int_{\Omega} \exp(\frac{1}{x^2+y^2+2})$ &$0.7508533738$ &$0.750772037320043$& $0.0000813364799567839$ \\ 
	\hline 
	$\int_{\Omega} \ln(\frac{1}{x^2+y^2+2})$ &-$0.4547712524$ &-$0.453623297839054$& $0.0011479545609460200$  \\ 
	\hline  \hline 
\end{tabular}
		\caption{\textbf{Comparison between numerical integration using $\delta_{-} ziti$ and the exact value, in the Cartesian case.}}
			\label{tab1}
\end{table}
\par 

\section{\textbf{Second strategy:  Polar coordinates.}} \label{poll}

In this section, we built all the necessary elements for the approximation $ \delta_{-}ziti $, using polar coordinates. The domain $\Omega=B(0,1)$ is represented using the polar coordinates , with the following parametrization:
$$ \forall (x,y)\in B(0,1), x=r\cos(\theta), y=r \sin(\theta), (r,\theta) \in [0,1] \times [0,2\pi].$$
\subsection{Construction of the method's tools using polar coordinates} \leavevmode \par 
This first part of the algorithm, compute $\lambda_{i}$ and $r_i$ in $[0,1].$\\
The construction's algorithm using the two variables $r$ and $\theta$  is given as follows:: \\
\begin{algorithm}
	\KwData \par  {~~~~~~~~~~~~\ \ \ $N_r=$ the root's number in the interval $[0,1]$ \\
		~~~~~~~~~~~~\ \ \ $h_r=\frac{1}{N_r-1}$ } \;
	\For {$i=1\cdots N_r$}{
		$x_{i}=a+(i-1)h_r$ \\
		$\varphi_{i}(r)=\Phi(\frac{r-x_i}{h_r})$
	} 
	$\alpha= \int_{0}^{1} \varphi_{1}(r) \varphi_{2}(r) dr $ \\
	$\beta= \int_{0}^{1} \varphi_{1}(r) \varphi_{1}(r) dr $ \\
	$\lambda_{1}=-\frac{\alpha}{\beta}$ \\
	\For {$i=1\cdots N_r-1$}{
		$\lambda_{i+1}=-\frac{\lambda_{1}}{2-\lambda_{1} \lambda_{i}}$
	}
	$\tilde{\Psi_1}(r)=\varphi_{1}(r)$ \\
	\For {$i=2\cdots N_r$}{
		$\tilde{\Psi_{i}}(r)=\varphi_{i}(r)+\lambda_{i-1} \tilde{\Psi}_{i-1}(r)$ \\
		
	} 
	\For {$i=1\cdots N_r$}{
		$\Psi_{i}(r)= \frac{\tilde{\Psi}_i(r)}{ \mid \mid \tilde{\Psi}_i (r)\mid \mid }$
	}
	\For {$i=1\cdots N_r$}{
		$\Lambda_i =\ln (-\lambda_{i})$, \\
		$P(y)=\Lambda_i y^4-2\Lambda_i y^3-\Lambda_iy^2+2(\Lambda_i-1)y+1=0$, \\
		$P(y^*_i)=0$, \\
		$r_i=x_i+h_r.y^*_i$
	}
\caption{Construction of $\delta_{-}ziti$'s tools, to compute the roots $r_i$}
\end{algorithm} \par 

\newpage
The second part of the algorithm permits to compute $\lambda_{i}^\theta$ and $s_i$ in $[0,2\pi].$
\begin{algorithm}
	\KwData \par  {~~~~~~~~~~~~\ \ \ $N_\theta=$ the root's number in the interval $[0,2\pi ]$.\\
		~~~~~~~~~~~~\ \ \ $h_\theta=\frac{2\pi}{N_\theta-1}$ } \;
	\For {$i=1\cdots N_\theta$}{
		$\theta_{i}=a+(i-1)h_\theta$ \\
		$\varphi_{i}(\theta)=\Phi(\frac{\theta-\theta_i}{h_\theta})$
	} 
	$\alpha^\theta= \int_{0}^{2 \pi } \varphi_{1}(\theta) \varphi_{2}(\theta) d\theta $ \\
	$\beta^= \int_{0}^{2 \pi } \varphi_{1}(\theta) \varphi_{1}(\theta) d\theta $ \\
	$\lambda_{1}^=-\frac{\alpha^\theta}{\beta^\theta}$ \\
	\For {$i=1\cdots N_\theta-1$}{
		$\lambda_{i+1}^=-\frac{\lambda_{1}^\theta}{2-\lambda_{1}^\theta \lambda_{i}^\theta}$
	}
	$\tilde{\Psi_1}(\theta)=\varphi_{1}(\theta)$ \\
	\For {$i=2\cdots N_\theta$}{
		$\tilde{\Psi_{i}}(\theta)=\varphi_{i}(\theta)+\lambda_{i-1}^\theta \tilde{\Psi}_{i-1}(\theta)$ \\
		
	} 
	\For {$i=1\cdots N_\theta$}{
		$\Psi_{i}(\theta)= \frac{\tilde{\Psi}_i(\theta)}{ \mid \mid \tilde{\Psi}_i (\theta)\mid \mid }$
	}
	\For {$i=1\cdots N_\theta$}{
		$\Lambda_i =\ln (-\lambda_{i}^\theta)$, \\
		$P(y^\theta)=\Lambda_i y^4-2\Lambda_i y^3-\Lambda_iy^2+2(\Lambda_i-1)y+1=0$, \\
		$P(y^{*^\theta}_i)=0$, \\
		$s_i=\theta_i+h_\theta.y^{*^\theta}_i$
	}
\caption{Construction of $\delta_{-}ziti$'s tools, to compute the roots $s_i$.}
\end{algorithm} \par 

\newpage

	The polar  set $(\Psi_{ij}(r,\theta))$  is defined by,
	\begin{equation}
	\Psi_{ij}(r,\theta):=\Psi_i(r).\Psi_j(\theta).
	\end{equation}

\subsection{Fundamental results: Numerical Integration} \leavevmode \par 
To test the previous strategy, we present in the following table, some numerical tests. We compare the exact value with the numerical one, using $\delta_{-} ziti$ method in the polar case.\par 

\renewcommand{\arraystretch}{2.5}
\begin{table}[htbp!]
	\centering
\begin{tabular}{ |c |  c |  c |    c  ||  }
	\hline \hline 
	$\displaystyle\frac{1}{2\pi} \int_{\Omega} f(x,y) dx dy$ & Ex\ \ \ & Ap\ \ \ & Error \\
	\hline
	$\displaystyle \int_{\Omega} (\frac{1}{x^2+y^2+1})^{\frac{1}{4}} dx dy$ &$0.4545285537$ \ \ \  &$0.454496459918650$ \ \ \  & $0.0000320937813496069$  \\ 
	\hline 
	$\displaystyle\int_{\Omega} \sqrt{(\frac{1}{x^2+y^2+1})} dx dy$ &$0.4142135624$ \ \ \  &$0.414440692467526$ \ \ \  & $0.0002271300675258380$  \\ 
	\hline 
	$\displaystyle\int_{\Omega} \exp(\frac{1}{x^2+y^2+2})dx dy$ &$0.7508533738$ \ \ \  &$0.750772037320043$\ \ \  & $0.0000813364799567839$ \\ 
	\hline 
	$ \displaystyle\int_{\Omega} \ln(\frac{1}{x^2+y^2+2})dx dy$ &$0.4547712524$ \ \ \  &-$0.453623297839054$\ \ \  & $0.0011479545609460200$  \\ 
	\hline 
	$\displaystyle\int_{\Omega}  x y dx dy$ & 0  &$0.000104893284924629$\ \ \  & $0.0001$\\
	\hline 
		$\displaystyle \int_{\Omega}  \frac{\ln(\sqrt{x^2+y^2})}{\sqrt{x^2+y^2}}dx dy$ & -1  &$-1.00011299531961$\ \ \  & $0.000112995319609510$\\
	\hline  \hline 
\end{tabular} 
		\caption{ \textbf{Comparison between numerical integration using $\delta_{-} ziti$ and the exact value, in the polar case.}}
			\label{tab2}
	\end{table}

Using the polar coordinates $r\in[0,1]$ and $\theta\in[0,2\pi]$ for some types of integrals, the following table shows us the error between exact and approximated solution founded using $\delta_{-}ziti$ method.
\newpage
\begin{table}[h]
	\centering
	\begin{tabular}{ |	c |  c  | c   |  c  ||  }
		\hline \hline 
			$\displaystyle \int_{\Omega} f(r,\theta) dr d\theta$ & Ex\ \ \ & Ap\ \ \ & Error \\
		\hline
	$\displaystyle \int_{\Omega} \frac{r.\sin(\theta)}{(r^2+1).t^{\frac{1}{3}}}$ &$0.2204366348$ \ \ \  &$.223358184762906$\ \ \  & $0.00292154976290612$  \\ 
\hline
$\displaystyle \int_{\Omega} \frac{\sin(\theta)}{r^{\frac{1}{2}}}$ &$0$ \ \ \  &$0.000796466856449210$\ \ \  & $0.000796466856449210$  \\ 
\hline
$\displaystyle \int_{\Omega} \frac{1}{\sqrt{2\pi r}}$ &$\frac{1}{\pi}$ \ \ \  &$\frac{0.994955526251694}{\pi}$\ \ \  & $0.00160570586425771$  \\ 
\hline
$\displaystyle \int_{\Omega} \frac{\ln(r)}{\sqrt{2\pi}}$ &$\frac{-0.500}{2\pi}$ \ \ \  &$\frac{-0.511242350773147}{2\pi}$\ \ \  & $0.00357855139409662$  \\ 
	\hline  \hline 
\end{tabular} 
\caption{ \textbf{Generalised integrals expressed with polar variables }}
\label{tabb2}
\end{table}
In the previous table, we remark that even we choose a generalised integral, like the example $\int_{\Omega} \frac{\sin(\frac{\theta}{2})}{x^{\frac{1}{2}}}$, which is in fact an operation of Riemann integral, we found a good approximation using $\delta_{-} ziti$ roots. For the two last examples, $\displaystyle \int_{\Omega} \frac{1}{\sqrt{2\pi r}}$ and $\displaystyle \int_{\Omega} \frac{\ln(r)}{\sqrt{2\pi}}$, other approximation methods (e.g.\, Simpson, Trapeze..) didn't give any result, which is an important point for our construction.

\section{\textbf{Numerical applications.}} \label{EDP}
\subsection{Elliptic PDE case : Poisson problem} \leavevmode 
\subsubsection{\textbf{The Cartesian case}}  \label{some} \leavevmode  \par 
In this section, let consider a Partial Differential Equation, which admits an exact solution and  we will compare it with the numerical one, using our method in the Cartesian case. Let $\Omega=B(0,1)$. The problem studied is given by:

\begin{subnumcases}{}
-\Delta u\ \ =f \ \ \ \ in \ \ \Omega, \label{kl} \\
	u(x,y)=0 \ \ \ \ in  \ \  \partial \Omega,    
	\label{MM}
\end{subnumcases}

with a given analytical solution $u_{ex}=1-x^2-y^2$, and the source term function is defined by  $f=4 $. \\
\subsubsection*{\textbf{The strong discretization}} \leavevmode \par 

The first step to approach the previous problem, is to multiply the  equation (\ref{kl}) by a test function $\Psi_{ij}$ and  to integrate the result over the domain $\Omega=B(0,1)$, which gives,
\begin{equation}
-\int_{\Omega} \Delta u(x,y).\Psi_{ij}(x,y) dx dy= \int_{\Omega} f(x,y).\Psi_{ij}(x,y) dx dy.
\end{equation}
Using  the theorem (\ref{integ}), we obtain the following scheme:

\begin{equation}
- \frac{\Delta u(r_i^j,s_j^i)}{\Psi_{ij}(r_i^j,s_j^i)}=\frac{f(r_i^j,s_j^i)}{\Psi_{ij}(r_i^j,s_j^i)}, \forall (r_i^j,s_j^i)\in \Omega.
\end{equation}
The next step, consists to approach the second derivative, which gives us the following scheme:

\begin{equation}
\left\lbrace
\begin{aligned}
&\frac{u_{i-1,j}+2u_{ij}+u_{i+1,j}}{(r_{i+1}^j-r_i^j)(r_i^j-r_{i-1}^j)}+\frac{u_{i,j-1}+2u_{ij}+u_{i,j+1}}{(s_{j+1}^i-s_j^i)(s_j^i-s_{j-1}^i)}=f(r_i^j,s_j^i),\ \ i,j=2\cdots N-1,\\
&u_{1,j}=u_{M,j}=0, \ \  j=1\cdots N,\\
&u_{i,1}=u_{i,M}=0, \ \  i=1\cdots N,\\
\end{aligned}
\right.
\label{Psgh}
\end{equation}
where, $N$ is the nods number in every internal segment (horizontally and vertically). At the end, we will have a global matrix, with $(N-2)\times(N-2)$ lines and $(N-2)\times(N-2)$ columns, defined as follows::
\[
M=
\left[
\begin{array}{c| c| c| c| c}
D^2 & A^3& 0& \cdots & 0 \\
\hline
A^2 & D^3& A^4& \cdots & 0\\
\hline
& & \vdots&  & \\
\hline
0& \cdots&  A^{n-3} & D^{n-2}& A^{n-1}\\
\hline 
0 & 0  & \cdots  &  A^{n-2} & D^{n-1}
\end{array}
\right]
\]
where $D^i$ is a $(N-2)\times (N-2)$ tri-diagonal matrix, defined by:
\begin{align*}
& D^{i+1}_{k,k} =-\frac{2\Psi_{i,k+1}}{dx_k.dx_{k-1}}-\frac{2\Psi_{i,k+1}}{dy_k.dy_{k-1}}, & i,k=1\cdots N-2, \\
& D^{i+1}_{k,k+1} =-\frac{2\Psi_{i,k+2}}{dy_k.dy_{k-1}}, & i,k=1\cdots N-2, \\
& D^{i+1}_{k-1,k}= -\frac{2\Psi_{i,k}}{dy_k.dy_{k-1}}, & i,k=1\cdots N-2, 
\end{align*}

and $A^i$ is a $(N-2)\times (N-2)$ diagonal matrix defined by:
\begin{equation*}
A^i_{k,k}=-\frac{\Psi_{i,k}}{dx_k.dx_{k-1}}, \ \  i=3\cdots N-1, \ \  and \ \ k=2\cdots N-1,
\end{equation*}

where, 
\begin{align*}
dx_k=r^j_k-r^j_{k-1}\raisepunct{,}\\
dy_k=r^i_k-r^i_{k-1}\raisepunct{,}
\end{align*}
therefore, we should resolve a simple system in the form $MX=F$, when $M$ is the global matrix defined previously,  $X$ in the unknown vector of size $(N-2)\times (N-2) $ and  $F$ is the source term vector of size $(N-2)\times (N-2) $.

\begin{rem}\leavevmode \par
To complete the resolution of the previous system, we must add boundary conditions ( homogeneous Dirichlet in this case).
\end{rem}
\subsubsection{\textbf{Numerical results}} \leavevmode \par 

Let $\Omega=B(0,1)$. In this case, we fix the points number in every single segment and  we vary the subdivision step. It is clear that the minimum of all the steps is obtained at the first segment (horizontally or vertically) and  the maximum is on the segment confused with the diameter of the disk, (\textit{i.e.}\, for two different intervals, horizontally or vertically, the associated step is not the same.)  \vskip 6pt

We are interested by the shape of the approximated solution with $\delta_{-} ziti$ scheme, using the Cartesian coordinates and  the segments approach.\\
For a fixed node's number in every segment (horizontal or vertical), $N=100$,
the numerical implementation of the scheme gives us an approximated solution, which is near to the exact one, given by $u_{ex}(x,y)=1-x^2-y^2,$ and $f(x,y)=4$.\\
 \begin{figure}[htbp!]
	\includegraphics[scale=0.6]{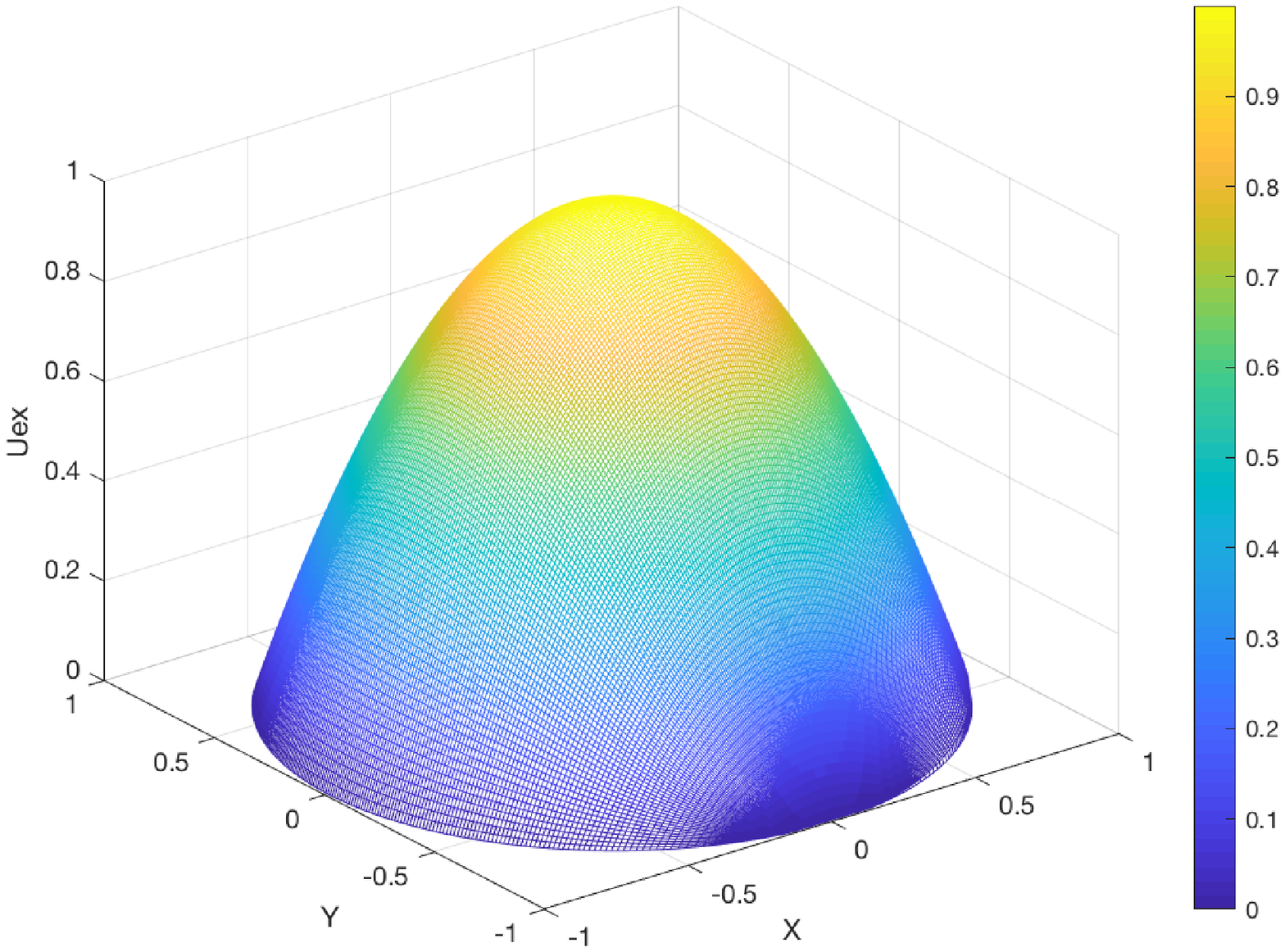}
	\caption{The exact solution}
	\label{ex_seg}
\end{figure}

\begin{figure}[htbp!]
	\includegraphics[scale=0.6]{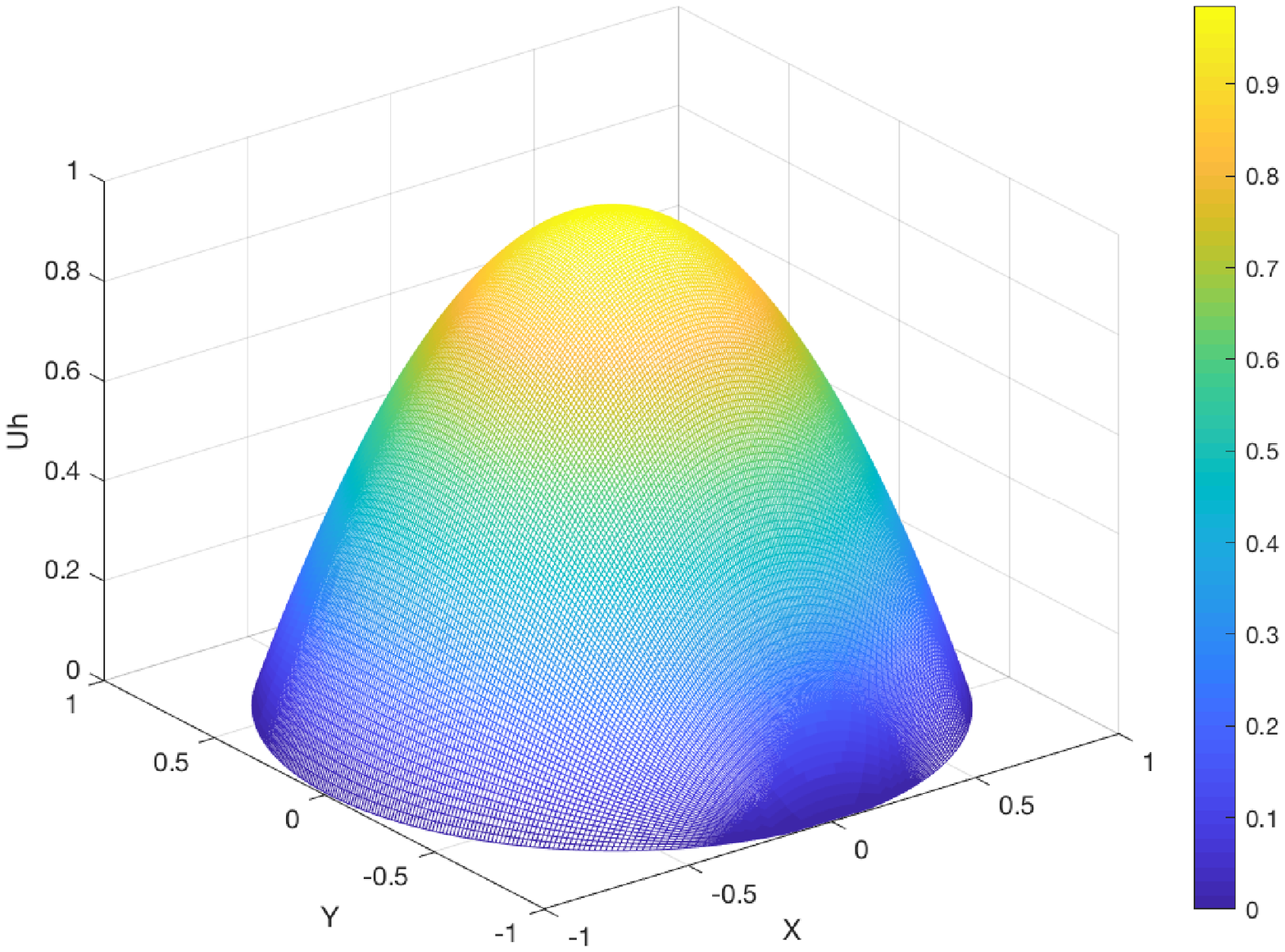}
	\caption{Approximated solution using $\delta_{-} ziti$}
	\label{ap_seg}
\end{figure}
\newpage
In the  table \ref{tab3}, we present the error between exact and approximated solution, with different values of nodes number $N$. \leavevmode \par 
\begin{table}[ht]
\begin{tabular}{ |c   |c|  c|     c |c ||    }
	\hline \hline 
	\centering
	$N$ \ \ \ \ \ \ \ \ \ \ & $h_{min}$\ \ \ \ \ \ \ \ \ \ & $h_{max}$\ \ \ \ \ \ \ \ \ \ & $Er_{max}$ & $Er_{mean}$ \\
	\hline
	$60$ \ \ \ \ \ \ \ \ \ \ & $0.008534$\ \ \ \ \ \ \ \ \ \ & $0.03333$\ \ \ \ \ \ \ \ \ \ & $0.01174$ & $0.00758$ \\
	\hline
	$100$ \ \ \ \ \ \ \ \ \ \ & $ 0.003979$\ \ \ \ \ \ \ \ \ \ & $0.02000$\ \ \ \ \ \ \ \ \ \ & $0.0046$& $0.00288$ \\
	\hline
	$150$ \ \ \ \ \ \ \ \ \ \ & $0.00217$\ \ \ \ \ \ \ \ \ \ & $0.01333$\ \ \ \ \ \ \ \ \ \ & $0.00204$& $0.0013169$ \\
	\hline
	$200$ \ \ \ \ \ \ \ \ \ \ & $0.001410$\ \ \ \ \ \ \ \ \ \ & $0.01000  $\ \ \ \ \ \ \ \ \ \ & $0.00167$& $0.0007514823$ \\
	\hline \hline 
\end{tabular} \par 
\caption{\textbf{The committed error using several values of nodes number $N$}}
\label{tab3}
\end{table} \par 

\newpage
where, 
\begin{align*}
Er_{max}& :=\underset{i,j}{\max}(|u_{ex}(i,j)-u(i,j))|,\\
Er_{mean}& :=\frac{1}{N} \sum_{i,j}|u_{ex}(i,j)-u(i,j))|.
\end{align*}

We present in the following subsection, a comparison between approximated solutions using Finite Element and $\delta_{-}ziti$ methods.
\subsubsection{Comparison with Finite elements method} \leavevmode \par 
Finite element method (FEM ) is a widely used analogy to resolve some types  of Partial Differential Equations. A large class of works was already done to resolve the Poisson problem using FEM, (see \cite{I} and \cite{J}). The starting point for the FEM is a PDE expressed in variational form. The basic recipe for turning a PDE into a variational problem is to multiply the equation by a test function $v$ and to integrate the resulting expression over all the domain $\Omega$: It is the common step between Galerkin analogy and $\delta_{-}ziti$. In this part, we present the approximated solution of the Poisson's problem defined in (\ref{kl}), using Finite Element Method.

	\begin{figure}[h]
			\includegraphics[scale=0.6]{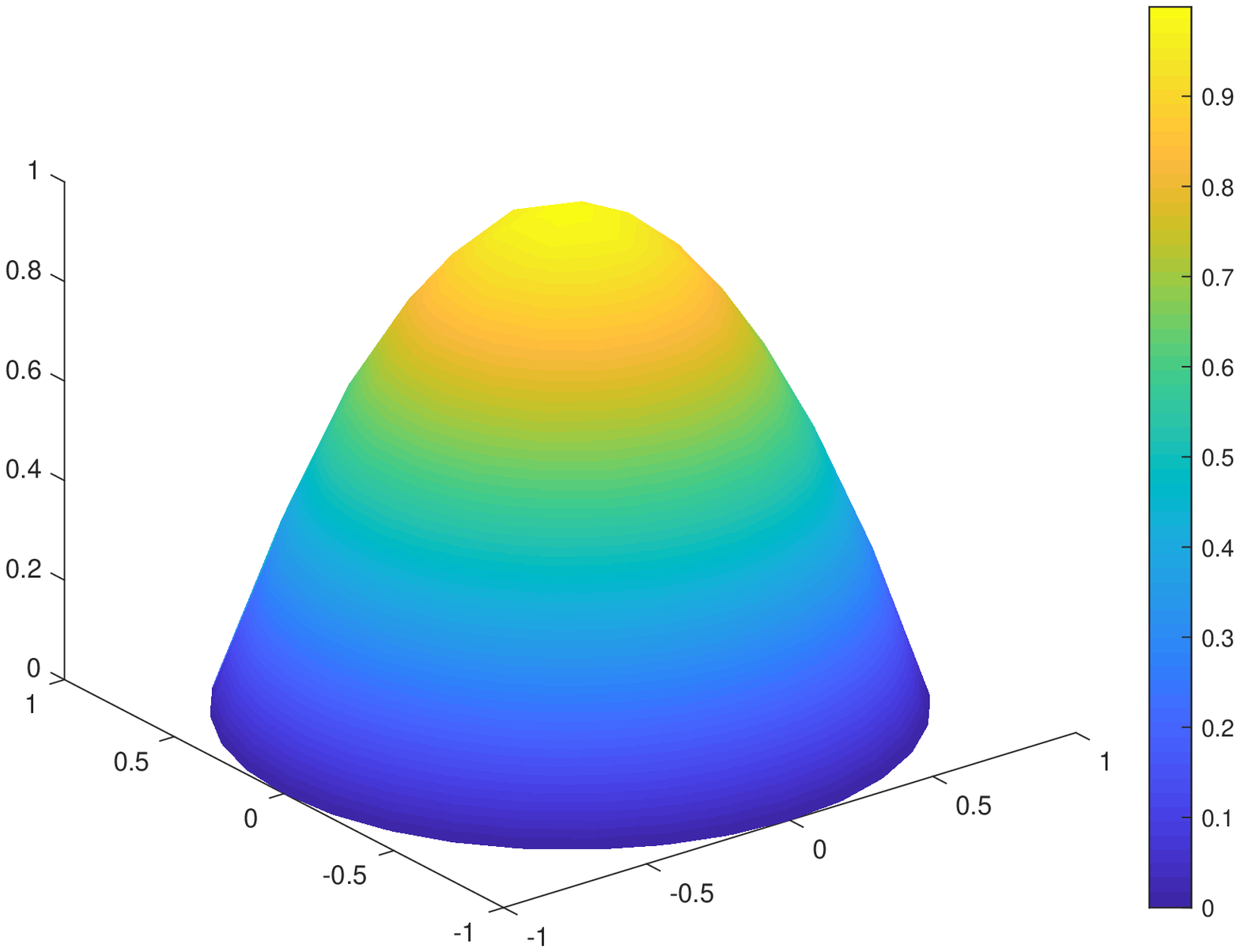}
			\caption{Approximated solution using finite element method }
	\end{figure}

\subsubsection{\textbf{The polar case}} \leavevmode \par 
Now, we consider the same partial differential equation defined before, which admits a polar analytical solution,  we will compare it  with the approximated one  founded using  $ \delta_{-}ziti$ method in the polar case.\par 
Let $ \Omega = B (0,1) $, the strategy presented consists in using the results of approximation in the mono-dimensional case and  taking into consideration the following function basis:
\begin{equation*}
\Psi_{ij}(r,\theta)=\Psi_i(r) \Psi_j(\theta),\ \  \forall i,j =1 \cdots N.
\end{equation*}

\begin{figure}[htbp!]
	\includegraphics[scale=0.4]{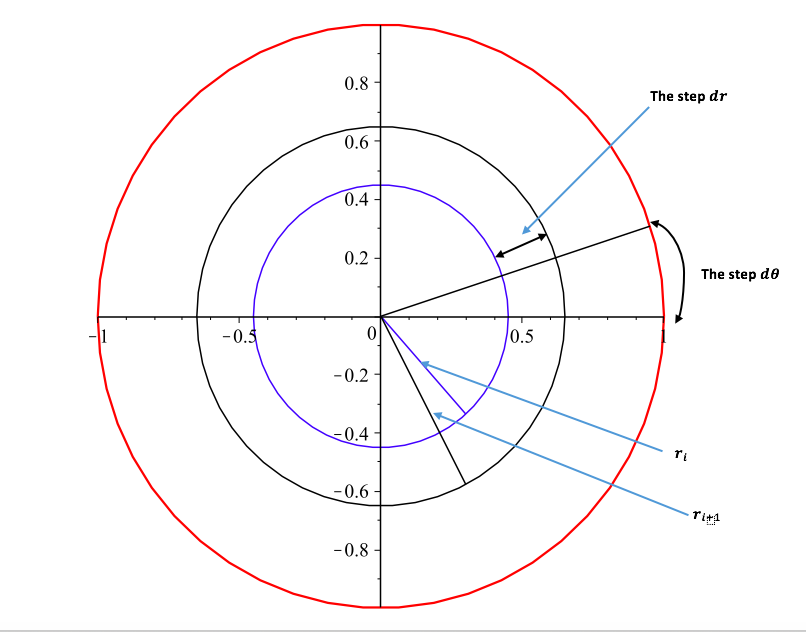}
	\caption{The domain $\Omega$ using polar parametrization}
	\label{fignh4}
\end{figure}
The problem presented in the previous subsection \ref{some} is equivalent of the polar one, expressed as follows:
\begin{subnumcases}{}
-\Delta u\ \ =\textbf{f} \ \ \ \ in \ \ \Omega, \label{klp} \\
u(r=1,\theta)=0 \ \ \ \forall \theta \in [0,2\pi],     
\label{key}
\end{subnumcases}

with 
\begin{align*}
& \Delta u := \frac{\partial^2 u}{\partial r^2} +\frac{1}{r}\frac{\partial u}{\partial r} +\frac{1}{r^2}\frac{\partial^2 u}{\partial \theta^2}, \\
& \textbf{f}(r,\theta)=f(r\cos\theta, r\sin\theta)=4
\label{pol}
\end{align*}
where $f$ id the source term defined in Cartesian problem, with an exact polar solution $u_{ex}=1-r^2$. Note that, the roots of the basic functions $ (\Psi_i (r)) $ will be noted $ r_i $ and  $ \theta_j $ are those associated with $ (\Psi_j (\theta)) $. \vskip 4pt
To obtain a numerical scheme using $\delta_{-} ziti$ method, we should multiply the equation \ref{klp} by a test function $\Psi_ {ij} (r, \theta) $ and  after, we use the strong result of approximation \ref{integ}, which gives:
\begin{equation}
\left\lbrace
\begin{aligned}
&\Delta_{ij}u=\textbf{f}(r_i,\theta_{j}),\ \ & i,j=2\cdots N-1,\\
&u_{1,j}=u_{2,j},\ \ & j=1\cdots N,\\
&u_{M,j}=0,\ \ & j=1\cdots N,\\
&u_{i,1}=u_{i,2},\ \ & i=1\cdots N,\\ 
&u_{i,M}=u_{i,M-1},\ \ & i=1\cdots N,\\ 
\end{aligned}
\right.
\label{Pss}
\end{equation}
with, $\Delta_{ij}u=\frac{u_{i-1,j}-2u_{ij}+u_{i+1,j}}{(r_{i+1}-r_i)(r_i-r_{i-1})}+\frac{1}{r_i} \frac{u_{i+1,j}-u_{ij}}{r_{i+1}-r_i} +\frac{1}{r_i^2} \frac{u_{i,j-1}+2u_{ij}+u_{i,j+1}}{(\theta_{j+1}-\theta_j)(\theta_j-\theta_{j-1})}$.\\ 
Therefore, the goal is to find $u(r,\theta)$ solution of polar problem given in (\ref{klp}) and (\ref{key}).\vskip 6pt
Like the Cartesian analogy, we resolve in this case, a simple system in the form $MX=F$, when $M$ is the global matrix defined previously, we should just ad the polar terms $\frac{1}{r}$ and $\frac{1}{r^2}$ in the corresponding terms of the matrices $A^i$ and $D^i$,  $X$ in the polar unknown vector of size $(N-2)\times (N-2) $ and  $F$ is the source term vector of size $(N-2)\times (N-2) $.\vskip 6pt
\textbf{Numerical tests.} \vskip 6pt 
For this test, we also consider the Poisson problem, but in the polar case.  Let $N_r=100$ be the root's number for the radius and  $N_\theta=100$ the other one for the angles. Two solutions, exact and approximated, are given by the following figures:
\begin{figure}[h]
	\includegraphics[scale=0.6]{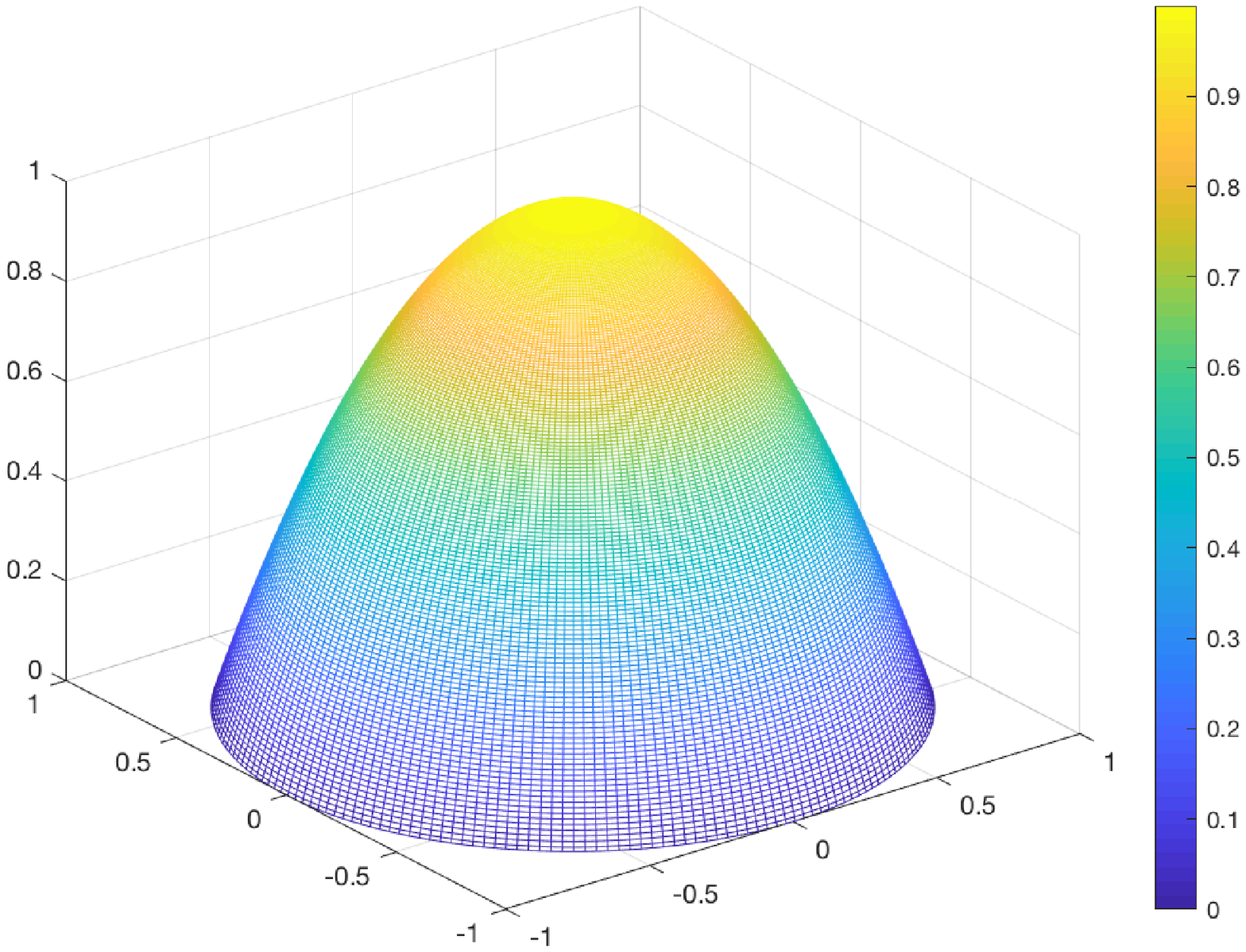}
	\caption{The exact solution of the problem}
	\label{fig3l}
\end{figure}
\begin{figure}[H]
	\includegraphics[scale=0.6]{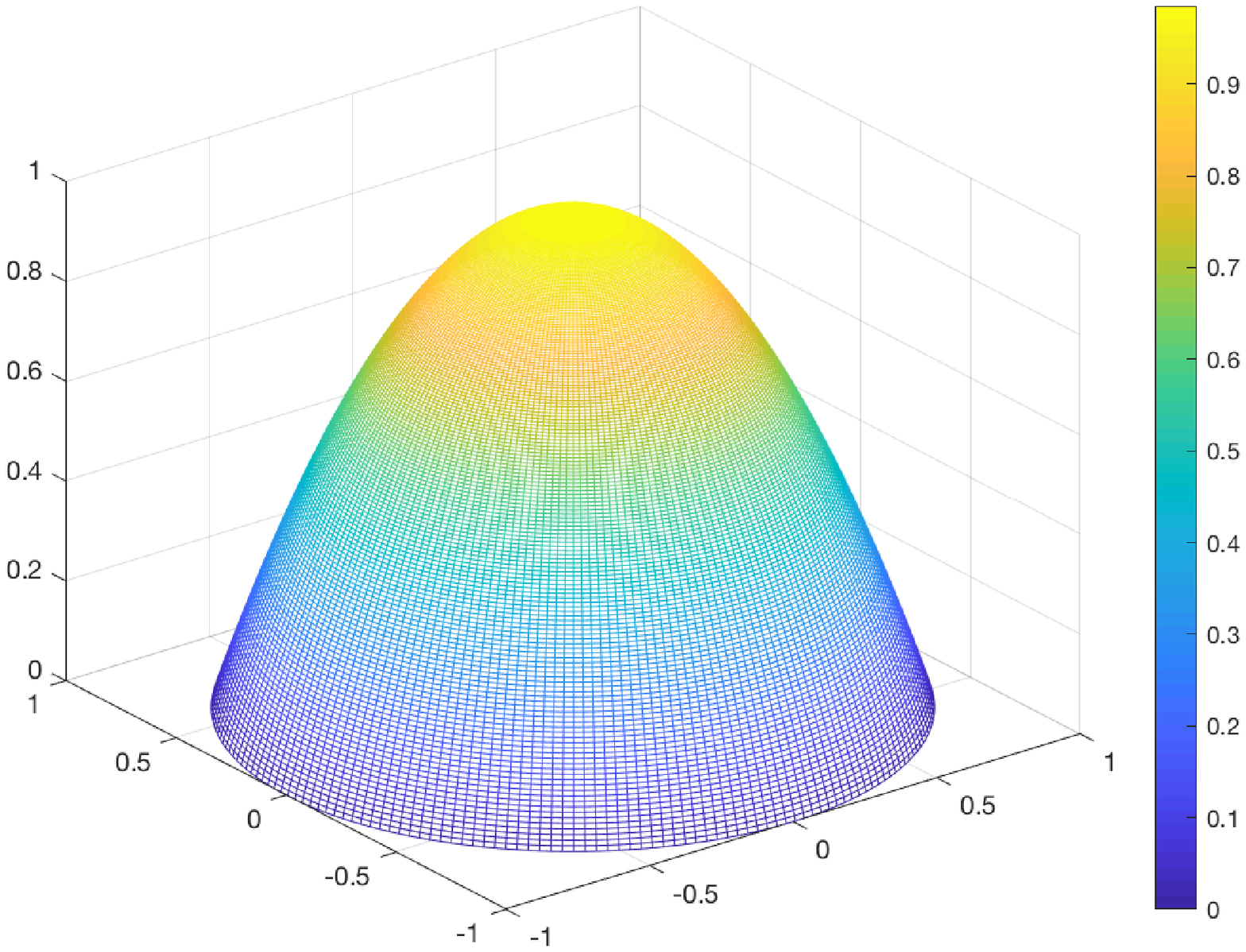}
	\caption{The approximated polar solution}
	\label{fig4}
\end{figure}
\subsubsection{Numerical Error}\leavevmode \par 
The  table \ref{tab4}, shows us the error between exact and approximated solution, using different $N_r$ and $N_\theta$. 
\begin{table}[h]
\begin{tabular}{| c  | c|  c |    c | c   ||}
	\hline \hline 
	\centering
	$N_r$ \ \ \ \ \ \ \ \ \ \ & $N_\theta$ \ \ \ \ \ \ \ \ \ \ & $h_r$\ \ \ \ \ \ \ \ \ \ & $h_\theta$\ \ \ \ \ \ \ \ \ \ & Error \\
	\hline
	$60$ \ \ \ \ \ \ \ \ \ \ & $60$ \ \ \ \ \ \ \ \ \ \ & $0.01666$\ \ \ \ \ \ \ \ \ \ & $0.0333 \pi $\ \ \ \ \ \ \ \ \ \ & $8.86.10^{-4}$ \\
	\hline
	$100$ \ \ \ \ \ \ \ \ \ \ & $100$ \ \ \ \ \ \ \ \ \ \ & $0.01$\ \ \ \ \ \ \ \ \ \ & $0.02 \pi $\ \ \ \ \ \ \ \ \ \ & $3.82.10^{-4}$ \\
	\hline
	$150$ \ \ \ \ \ \ \ \ \ \ & $150$ \ \ \ \ \ \ \ \ \ \ & $0.00666$\ \ \ \ \ \ \ \ \ \ & $0.0133 \pi $\ \ \ \ \ \ \ \ \ \ & $1.91.10^{-4}$ \\
	\hline
	$200$ \ \ \ \ \ \ \ \ \ \ & $200$ \ \ \ \ \ \ \ \ \ \ & $0.005$\ \ \ \ \ \ \ \ \ \ & $0.01 \pi $\ \ \ \ \ \ \ \ \ \ & $1.15.10^{-4}$ \\
	\hline \hline 
\end{tabular} \par 
\caption{\textbf{The infinite error, using several values of $N_r$ and $N_\theta$ }}
	\label{tab4}
\end{table}

\begin{figure}[h] 
	\begin{minipage}[t]{6cm} 
		\includegraphics[width=6cm]{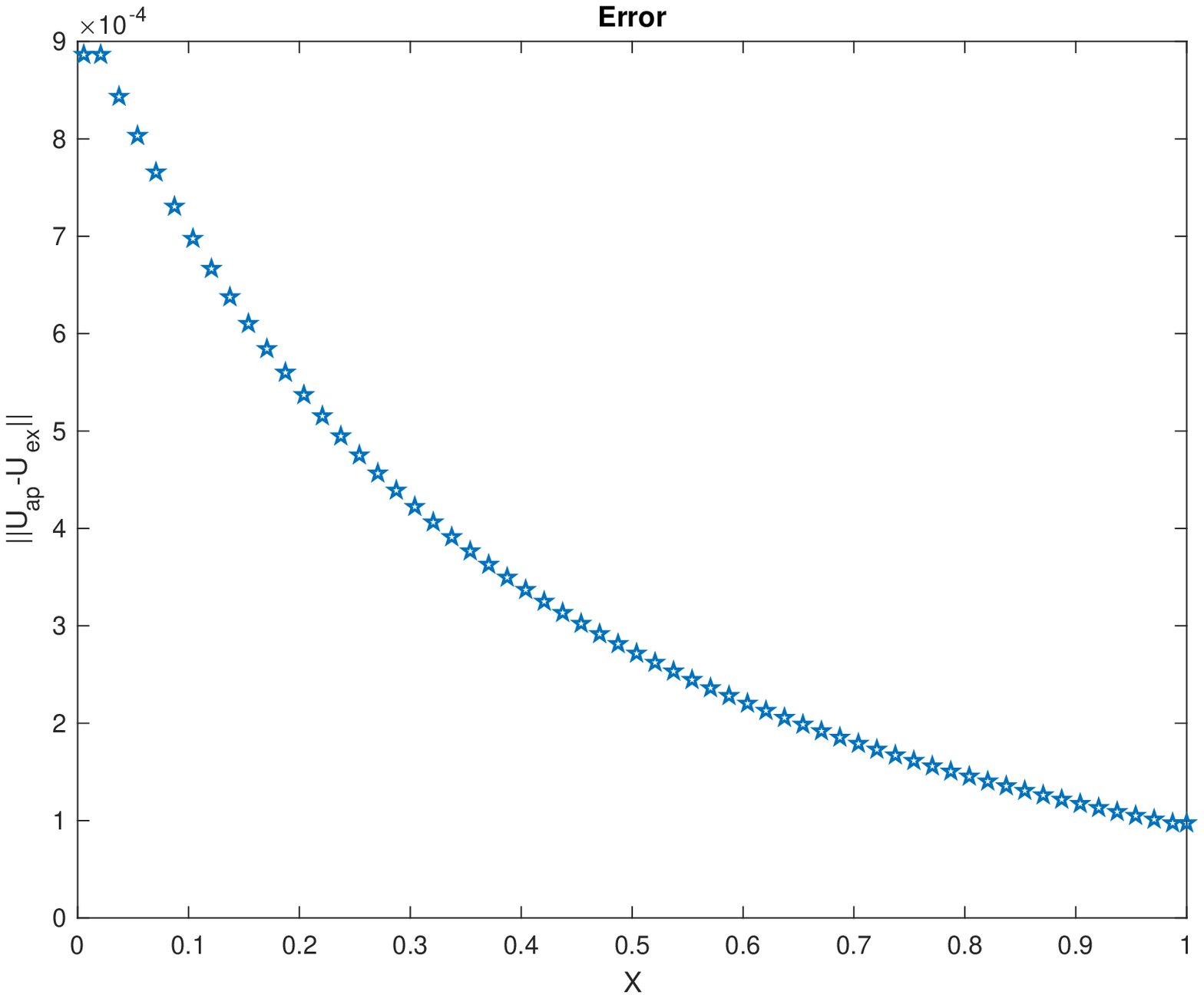}
		\title{The error for $N_r=60$.} 
	\end{minipage} 
	\begin{minipage}[t]{6cm}
		\centering \includegraphics[width=6cm]{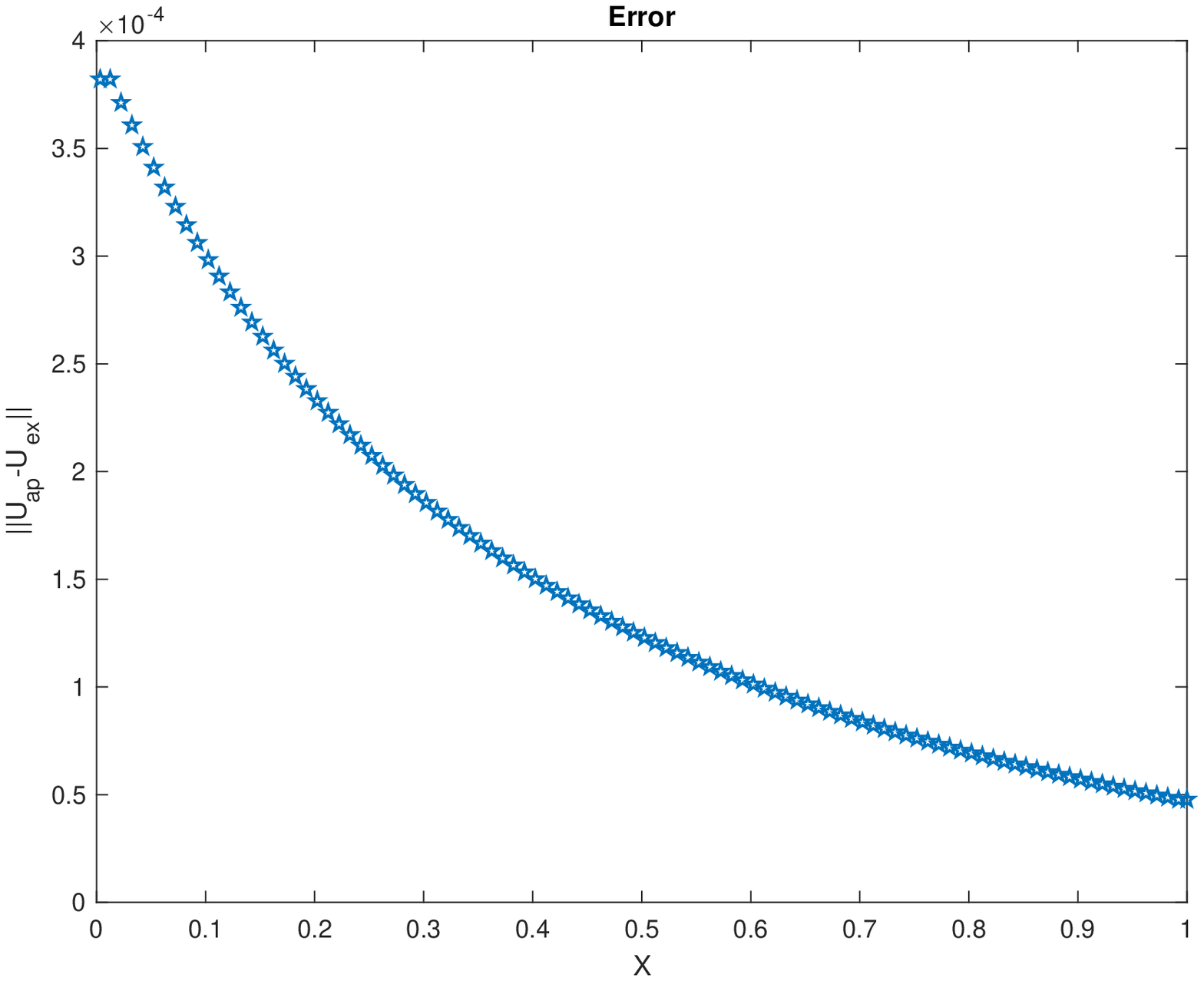}
			\title{The error for $N_r=100$.} 
	\end{minipage} 
	\begin{minipage}[t]{6cm} \centering
		\includegraphics[width=6cm]{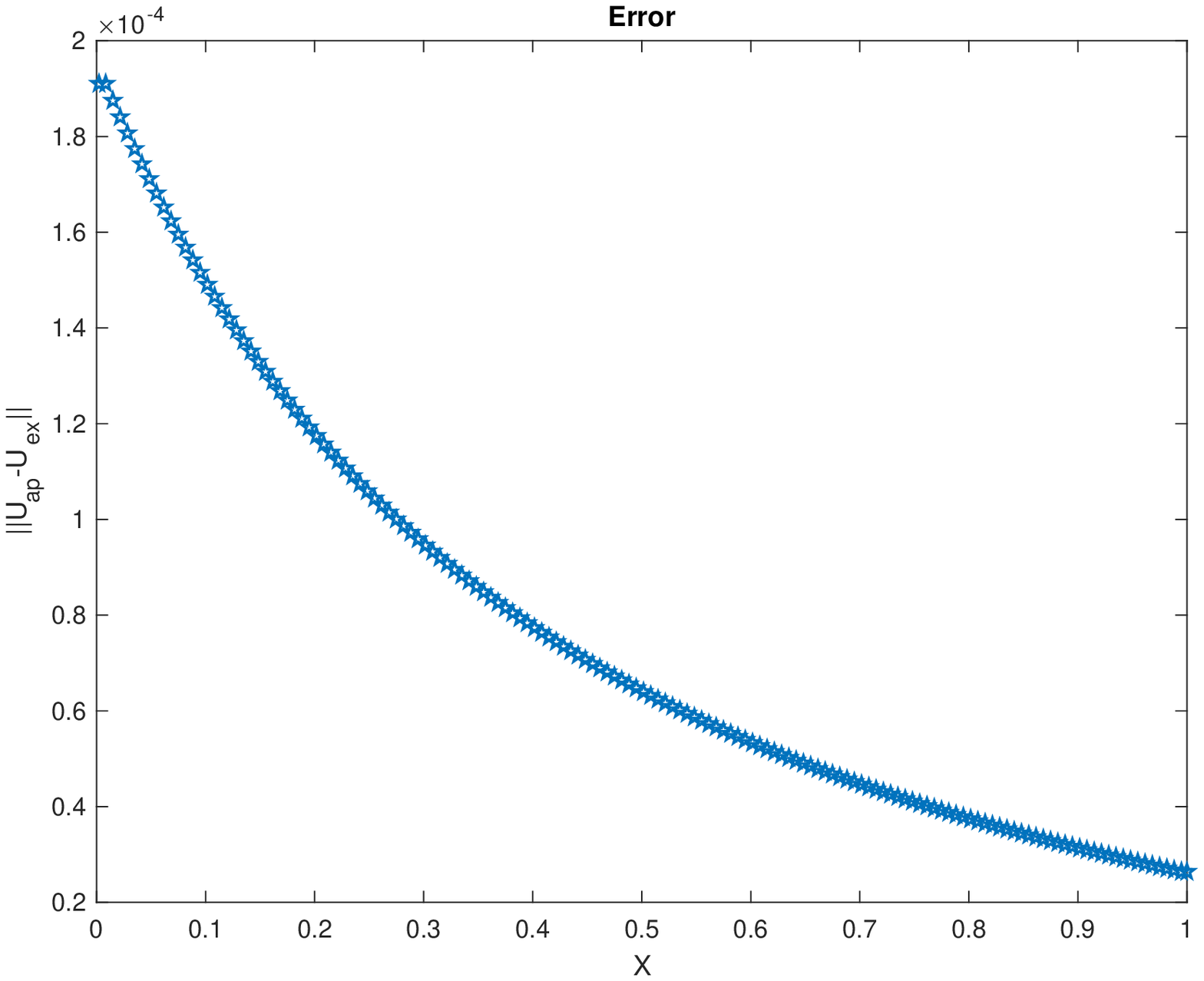}
		\title{The error for $N_r=150$.} 
	\end{minipage} 
	\begin{minipage}[t]{6cm}
		\centering \includegraphics[width=6cm]{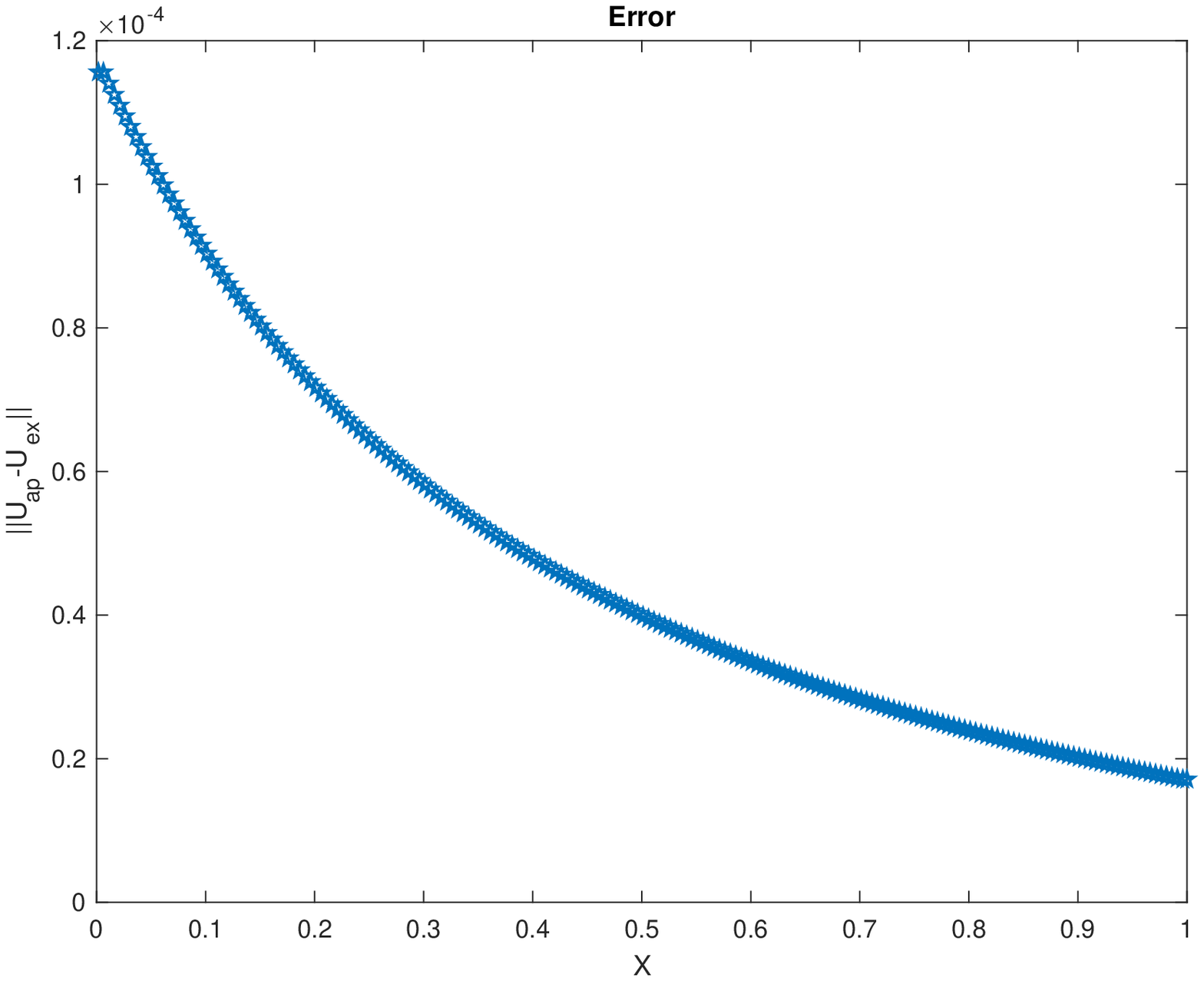}
			\title{The error for $N_r=200$.} 
	\end{minipage} 
\end{figure}

\newpage
\subsection{Parabolic PDE case: Heat equation} \leavevmode \par 
This section is devoted to the application of  $ \delta_{-} ziti $ method, on a diffusion equation, in a disk domain. The heat equation describes the   distribution of heat   (or variation in temperature)  in a  given region over time.  For  a  function, $u(t,x,y)$ (respectively $u(t,r,\theta)$) of  two  spatial  variables  $(x,y)$ in  the Cartesian case ($(r,\theta)$ in the polar case) and the time variable t, the heat equation is given by: 
\begin{subnumcases}{}
\frac{\partial u}{\partial t}  +D \Delta u=f, \ \ in\ \  \Omega, \label{kkl}\\
u(t,x,y)=0,\ \ (x,y)\in \partial \Omega,\\ 
u(0,x,y)=u_0(x,y) \geq 0.
\end{subnumcases}

Using the same analogy applied in the previous sections, we multiply the equation  \ref{kkl} by a test function $\Psi_{ij}$, after we integrate over the domain $\Omega$. It remains just the direct application of our approximations formulas given in \ref{integ}. The numerical scheme in the polar case is presented as follows:

\begin{equation*}
\left\lbrace
\begin{aligned}
&u_{ij}^{n+1}=u_{ij}^n-dt.D.\Delta_{ij}^nu+ dt.f_{i,j}^n,\ \ i,j=2\cdots N-1,\\
&u_{1,j}^{n+1}=u_{2,j}^{n+1},  \ \ \ \ \ \  \ \ \ \ \ \ j=1 \cdots N,\\
&u_{M,j}^{n+1}=0, \ \ \ \ \ \  \ \ \ \ \ \ \ \ \ \ j=1 \cdots N, \\
&u_{i,1}^{n+1}=u_{i,2}^{n+1},  \ \ \ \ \ \  \ \ \ \ \ \ i=1 \cdots N,\\
&u_{i,M}^{n+1}=u_{i,M-1}^{n+1}, \ \ \ \ \ \ \ \ \   i=1 \cdots N, \\
\end{aligned}
\right.
\label{sch_pol}
\end{equation*}
where, $\Delta_{ij}^n=\frac{u_{i-1,j}^n-2u_{ij}^n+u_{i+1,j}^n}{(r_{i+1}-r_i)(r_i-r_{i-1})}+\frac{u_{i,j-1}^,+2u_{ij}^n+u_{i,j+1}^n}{(\theta_{j+1}-\theta_j)(\theta_j-\theta_{j-1})}.$\\ \\
 Our goal is to find an approximated solution, near of the exact one, which verify the boundary conditions. The following function
\begin{equation*}
u_{ex}(t,x,y)=(1-x^2-y^2)\exp(t),
\end{equation*}
is an exact solution of the heat equation, with, 

\begin{equation*}
D=1, \ \ \ and \ \ f(t,x,y)=(-3-x^2-y^2))\exp(t).
\end{equation*}

The following figures, shows as the allure of exact and approximated solution, at a given finite time. Note that, with a simple variable changing, we can use a very high $t_f$, which is very useful to reduce the number of time iterations.
\begin{figure}[h]
	\includegraphics[scale=0.6]{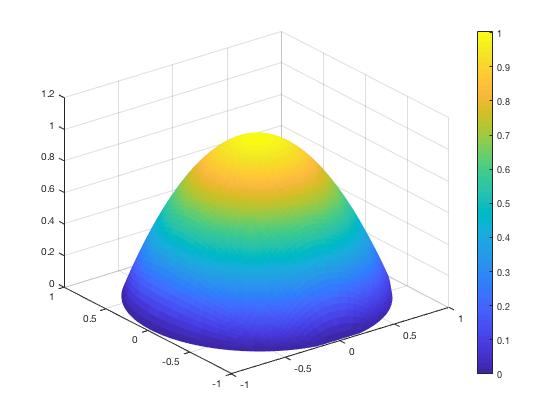}
	\caption{The exact solution }
	\label{fig5}
\end{figure}
\begin{figure}[h]
	\includegraphics[scale=0.6]{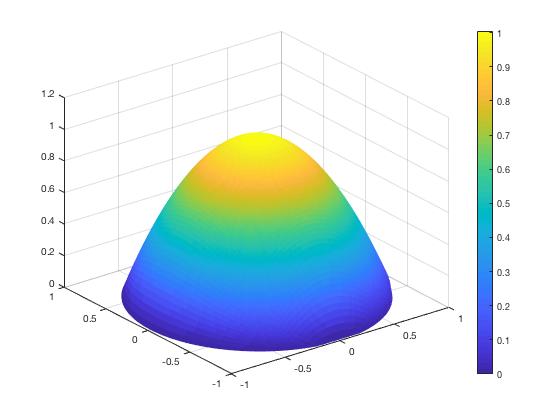}
	\caption{The approximated  solution using $\delta_{-}ziti$ method }
	\label{fig6}
\end{figure}
\newpage
In the  table \ref{tab5}, we present the error between exact and approximated global  solution, at a given finite time.  In this test, we take the following parameters:\\
\begin{equation*}
t_f=1.63.10^{16}, \mu=0.1,
\end{equation*}
\newpage
\begin{table}[h]
\begin{tabular}{ |c |  c  |c    | c |c ||    }
	\hline \hline 
	\centering
	$N$ \ \ \ \ \ \ \ \ \ \ & $h_{min}$\ \ \ \ \ \ \ \ \ \ & $h_{max}$\ \ \ \ \ \ \ \ \ \ & $Er_{max}$ & $Er_{mean}$ \\
	\hline
	$60$ \ \ \ \ \ \ \ \ \ \ & $0.008534$\ \ \ \ \ \ \ \ \ \ & $0.03333$\ \ \ \ \ \ \ \ \ \ & $2.22431.10^{-5}$ & $1.413039.10^{-5}$ \\
	\hline
	$100$ \ \ \ \ \ \ \ \ \ \ & $ 0.003979$\ \ \ \ \ \ \ \ \ \ & $0.02000$\ \ \ \ \ \ \ \ \ \ & $1.0174501 .10^{-6}$& $5.453807.10^{-7}$ \\
	\hline
	$150$ \ \ \ \ \ \ \ \ \ \ & $0.00217$\ \ \ \ \ \ \ \ \ \ & $0.01333$\ \ \ \ \ \ \ \ \ \ & $1.00359253.10^{-6}$& $5.36991485.10^{-7}$ \\
	\hline
	$200$ \ \ \ \ \ \ \ \ \ \ & $0.001410$\ \ \ \ \ \ \ \ \ \ & $0.01000  $\ \ \ \ \ \ \ \ \ \ & $9.9973690.10^{-7}$& $5.3444320.10^{-7}$ \\
	\hline \hline 
\end{tabular} \par 
\caption{\textbf{The committed error at a given finite time}}
	\label{tab5}
\end{table}

\subsubsection{Comparison with Finite Element Method} \leavevmode \par 

We present the approximated solution given by  the finite elements method. In this direction, a large body of works was already done, see \cite{I}. 
	\begin{figure}[htbp]
		\includegraphics[scale=0.6]{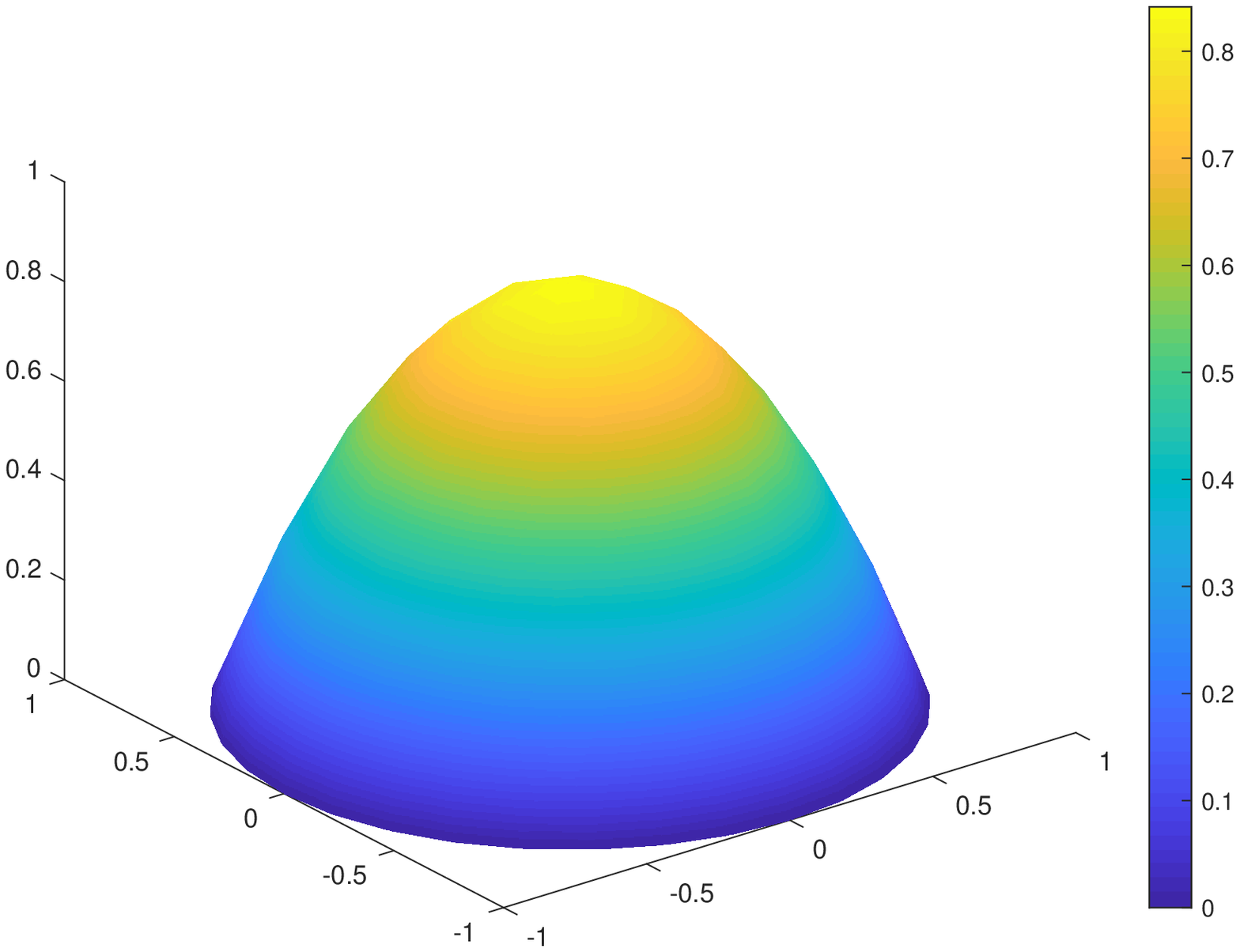}
		\caption{Approximated solution using finite elements method}
	\end{figure}

\section{\textbf{Conclusion}}
In the current work, we applied the $\delta_{-} ziti$ method in a disk domain, to calculate numerically some type of integrals,  to solve the Poisson and  Heat problem, using two strategies.  Firstly, concerning the approximated solution of a PDE, we start as the  Galerkin method, by constructing a weak formulation of the problem, then we use the roots of our orthonormal basis functions. Since this last  goes to Dirac function ( in the distribution sense), we can say that the $\delta_{-} ziti$ is a mix between Galerkin and the particular method. As a conclusion, $\delta_{-} ziti$ permits us to use two strategies in the case $\Omega=B(0,1)$, injecting the work already done in the mono-dimensional case (see \cite{Aa} and \cite{A}). The result is impressive, in fact:\vskip 6pt
\ \ \  $\bullet$ $\delta_{-} ziti$ is a fast scheme, precise, and gives an admissible solution.\vskip 3 pt
 \ \ \  $\bullet$ In the case of the Heat equation, the CFL condition of stability is near of $0.9$, and the numerical solution exist globally ($t_f$ goes to $10^{16}$). \\

\end{document}